\documentclass{conm-p-l}

\newtheorem{theorem}{Theorem}[section]
\newtheorem{lemma}[theorem]{Lemma}

\theoremstyle{definition}
\newtheorem{definition}[theorem]{Definition}

\theoremstyle{remark}

\numberwithin{equation}{section}

\newcommand{\E}{{\mathcal E}}

\newcommand{\e}{\epsilon}

\newcommand{\N}{{\mathcal N}}

\newcommand{\tQ}{\tilde{Q}}
\newcommand{\vQ}{\vec{Q}}

\newcommand{\C}{{\mathcal C}}
\renewcommand{\k}{\kappa}
\newcommand{\ga}{\gamma}

\newcommand{\hS}{\hat{S}}

\newcommand{\Dl}{\Delta}
\renewcommand{\th}{\theta}
\newcommand{\ra}{\rightarrow}
\newcommand{\al}{\alpha}
\newcommand{\be}{\beta}
\newcommand{\sg}{\sigma}
\newcommand{\Sg}{\Sigma}

\newcommand{\pa}{\partial}
\newcommand{\z}{\zeta}
\newcommand{\hQ}{\hat{Q}}

\newcommand{\hx}{\hat{x}}
\newcommand{\hz}{\hat{z}}

\newcommand{\La}{\Lambda}

\newcommand{\la}{\lambda}

\newcommand{\nid}{\noindent}

\newcommand{\om}{\omega}
\newcommand{\Om}{\Omega}

\newcommand{\tx}{\tilde{x}}
\newcommand{\ty}{\tilde{y}}
\newcommand{\tz}{\tilde{z}}
\newcommand{\vtQ}{\vec{\tilde{Q}}}

\newcommand{\W}{{\mathcal W}}
\renewcommand{\L}{{\mathcal L}}
\renewcommand{\O}{{\mathcal O}}
\newcommand{\T}{{\mathcal T}}

\def\maprightu#1{\smash{
    \mathop{\longrightarrow}\limits^{#1}}}
\def\maprightd#1{\smash{
    \mathop{\longrightarrow}\limits_{#1}}}
\def\mapdownl#1{
    \llap{$\vcenter{\hbox{$\scriptstyle#1$}}$}\Big\downarrow}
\def\mapdownr#1{\Big\downarrow
    \rlap{$\vcenter{\hbox{$\scriptstyle#1$}}$}}

\begin{document}

\title{Existence of Chaos for a Singularly Perturbed NLS Equation}

\author{Yanguang (Charles)  Li}
\address{Department of Mathematics, University of Missouri, 
Columbia, MO 65211}
\curraddr{}
\email{cli@math.missouri.edu}
\thanks{}


\date{}



\begin{abstract}
The work \cite{Li99} is generalized to the singularly perturbed 
nonlinear Schr\"odinger (NLS) equation of which the regularly 
perturbed NLS studied in \cite{Li99} is a mollification. Specifically,
the existence of Smale horseshoes and Bernoulli shift dynamics is 
established in a neighborhood of a symmetric pair of Silnikov homoclinic 
orbits under certain generic conditions, and the existence of the 
symmetric pair of Silnikov homoclinic orbits has been proved in 
\cite{Li01}. The main difficulty in the current horseshoe construction
is introduced by the singular perturbation $\e \pa_x^2$ which turns 
the unperturbed reversible system into an irreversible system. It turns
out that the equivariant smooth linearization can still be achieved, 
and the Conley-Moser conditions can still be realized.
\end{abstract}

\maketitle








\section{Introduction}

Consider the singularly perturbed nonlinear Schr\"odinger (NLS) equation,
\begin{equation}
iq_t = q_{\z\z} +2 [ |q|^2 - \om^2] q +i \e [q_{\z\z} - \al q +\be ] \ ,
\label{pnls}
\end{equation}
where $q = q(t,\z)$ is a complex-valued function of the two real 
variables $t$ and $\z$, $t$ represents time, and $\z$ represents
space. $q(t,\z)$ is subject to periodic boundary condition of period 
$2 \pi$, and even constraint, i.e. 
\[
q(t,\z + 2 \pi) = q(t,\z)\ , \ \ q(t,-\z) = q(t,\z)\ .
\]
$\al >0$ and $\be >0$ are constants, and $\e > 0$ is the perturbation 
parameter. For simplicity of presentation, we restrict $\om$ by 
$\om \in (1/2, 1)$. In \cite{Li01}, the following theorem on the 
existence of Silnikov homoclinic orbits was proved.
\begin{theorem}
There exists a $\e_0 > 0$, such that for any $\e \in (0, \e_0)$, there 
exists a codimension 1 surface in the external parameter space 
$(\alpha,\beta, \om) \in \mathbb{R}^+\times  \mathbb{R}^+\times 
\mathbb{R}^+$ where $\om \in (\frac{1}{2}, 1)/S$, $S$ is a finite subset, and 
$\al \om < \be$. For any $(\alpha ,\beta, \omega)$ on the codimension 1
surface, the singularly perturbed nonlinear Schr\"odinger equation 
\eqref{pnls} possesses a symmetric pair of Silnikov homoclinic orbits 
asymptotic 
to a saddle $Q_\epsilon$. The codimension 1 surface has the approximate 
representation given by $\al = 1/\k(\om)$, where $\k(\om)$ is plotted 
in Figure \ref{kappa}.
\label{horbit}
\end{theorem}
\begin{figure}
\vspace{3.0in}
\caption{The graph of $\k(\om)$.}
\label{kappa}
\end{figure}
Notice that if $q(t,\z)$ is a homoclinic orbit, then $q(t,\z+\pi)$ is 
another homoclinic orbit. Thus $q(t,\z)$ and $q(t,\z+\pi)$ form a 
symmetric pair of homoclinic orbits. Based upon the above theorem, we 
will construct Smale horseshoes in a neighborhood of the symmetric pair 
of homoclinic orbits. The construction is a generalization of that in 
\cite{Li99} where the singular perturbation $\e \pa_x^2$ is mollified 
into a bounded Fourier multiplier. The main difficulty in the current 
horseshoe construction is introduced by the singular perturbation 
$\e \pa_x^2$ which turns the unperturbed reversible system into an 
irreversible system. Specifically, denote by $F^t_\e$ the evolution 
operator of the singularly perturbed nonlinear Schr\"odinger equation 
\eqref{pnls}. When $\e =0$, $F^t_0$ is a group. When $\e >0$, $F^t_\e$
is only a semigroup. It turns out that the equivariant smooth linearization 
can still be achieved, and the Conley-Moser conditions can still be realized.
Of course, one has to replace the inverse of the evolution operator
$F^t_\e$ by preimage. The article is organized as follows: In section 2, 
we present equivariant smooth linearization. In section 3, we present 
the Poincar\'e map and its representation. In section 4, the fixed points of 
the Poincar\'e map is studied. In section 5, we present the existence 
of chaos. Finally, in section 6, numerical evidence for the generic 
conditions is presented.

\section{Equivariant Smooth Linearization}

The symmetric pair of Silnikov homoclinic orbits is asymptotic to the 
saddle $Q_\e =\sqrt{I} e^{i\th}$, where
\begin{equation}
I=\omega^2-\epsilon \frac{1}{2\omega}\sqrt{\beta^2-\alpha^2\omega^2}+\cdots ,
\quad \cos \theta  =\frac{\alpha \sqrt{I}}{\beta}, \quad \theta \in 
(0,\frac{\pi}{2}).
\label{Qec}
\end{equation}
Its eigenvalues are
\begin{equation}
\la_n^\pm = -\e [\al +n^2]\pm 2 \sqrt{(\frac{n^2}{2} +\om^2-I)(3I -\om^2 -
\frac{n^2}{2} )}\ , 
\label{Qev}
\end{equation}
where $n=0,1,2, \cdots $, $\om \in (\frac{1}{2}, 1)$, and $I$ is given in 
(\ref{Qec}). The crucial points to notice are: (1). only $\la_0^+$ and 
$\la_1^+$ have positive real parts, $\mbox{Re}\{ \la_0^+\} < \mbox{Re}\{ 
\la_1^+\} $; (2). all the other eigenvalues have negative real parts among 
which the absolute value of $\mbox{Re}\{ \la_2^+\}=\mbox{Re}\{ \la_2^-\}$ 
is the smallest; (3). $|\mbox{Re}\{ \la_2^+\}| < \mbox{Re}\{ \la_0^+\}$. 
Actually, items (2) and (3) are the main characteristics of Silnikov 
homoclinic orbits.
\begin{lemma}
For any fixed $\e \in (0,\e_0)$, let $E_\e$ be the codimension 1 surface in 
the external parameter space, on which the symmetric pair of Silnikov 
homoclinic orbits are supported (cf: Theorem \ref{horbit}). For almost 
every $(\al,\be,\om)\in E_\e$, the eigenvalues $\la_n^\pm$ (\ref{Qev}) 
satisfy the nonresonance condition of Siegel type: There exists a natural 
number $s$ such that for any integer $n \geq 2$,
\[
\bigg | \La_n - \sum_{j=1}^r \La_{l_j} \bigg | \geq  1/r^s\ ,
\]
for all $r =2,3, \cdots, n$ and all $l_1, l_2, \cdots, l_r \in \mathbb{Z}$, 
where $\La_n = \la_n^+$ for $n \geq 0$, and $\La_n = \la_{-n-1}^-$ for $n <0$.
\end{lemma}
\nid
Proof. The same proof as in \cite{Li99} can be carried through here. QED

Thus, in a neighborhood of $Q_\e$, the singularly perturbed NLS (\ref{pnls}) 
is analytically equivalent to its linearization at $Q_\e$ \cite{Nik86}. In 
terms of eigenvector basis, (\ref{pnls}) can be rewritten as 
\begin{eqnarray}
\dot{x} &=& -a x - b y +\N_x(\vec{Q}), \nonumber \\
\dot{y} &=& b x - a y +\N_y(\vec{Q}), \nonumber \\
\dot{z}_1 &=& \ga_1 z_1 +\N_{z_1}(\vec{Q}), \label{nfeq} \\
\dot{z}_2 &=& \ga_2 z_2 +\N_{z_2}(\vec{Q}), \nonumber \\
\dot{Q} &=& LQ +\N_{Q}(\vec{Q}); \nonumber
\end{eqnarray}
where $a= -\mbox{Re}\{\la^+_2\}$, $b=\mbox{Im}\{\la^+_2\}$, 
$\ga_1 = \la^+_0$, $\ga_2 = \la^+_1$; $\N$'s vanish identically 
in a neighborhood $\Om$ of $\vec{Q} =0$, $\vec{Q} = (x,y,z_1,z_2,Q)$, $Q$ 
is associated with the rest of eigenvalues, $L$ is given as
\[
LQ = -i Q_{\z\z} -2 i [(2|Q_\e|^2-\om^2)Q+Q^2_\e\bar{Q}]+
\e [-\al Q+Q_{\z\z}]\ ,
\]
and $Q_\e$ is given in (\ref{Qec}). The following theorem on well-posedness 
is standard \cite{Kat75}. Let $F^t$ ($0 \leq t < \infty$) be the 
evolution operator of the singularly perturbed NLS (\ref{nfeq}), and $H^s$ 
be the Sobolev space. 
\begin{theorem}
For any $s \geq 1$, and any $\vQ_0 \in H^{s+2}$, $F^t(\vQ_0) \in 
C^0([0,\infty);H^{s+2}) \cap C^1([0,\infty);H^{s})$. For any fixed 
$t \in [0,\infty)$, $F^t$ is a $C^2$ map in $H^s$.
\label{WPT}
\end{theorem}

\section{The Poincar\'e Map and Its Representation}

Denote by $h_k$ ($k=1,2$) the symmetric pair of Silnikov homoclinic orbits. 
The symmetry $\sg$ of half spatial period shifting has the new 
representation in terms of the new coordinates
\begin{equation}
\sg \circ (x,y,z_1,z_2,Q) = (x,y,z_1,-z_2,\sg \circ Q).
\label{symm}
\end{equation}
We have the following facts about the homoclinic orbits:
\begin{enumerate}
\item The homoclinic orbits are classical solutions,
\item As $t \ra -\infty$, the homoclinic orbits are tangent to the positive 
$z_1$-axis at $\vQ =0$.
\end{enumerate}
The same proof as in \cite{Li99} works here for item (2). Since $a$ is the 
smallest attracting rate, we assume that
\begin{itemize}
\item (A1). As $t \ra +\infty$, the homoclinic orbits are tangent to the 
($x,y$)-plane at $\vQ = 0$.
\end{itemize}
The Poincar\'e section is defined as in \cite{Li99}.
\begin{definition}
The Poincar\'e section $\Sg_0$ is defined by the constraints:
\begin{eqnarray*}
& & y=0,\ \eta \exp \{ -2 \pi a/b \} < x < \eta; \\
& & 0 < z_1 <\eta,\ -\eta < z_2 <\eta, \ \|Q\| < \eta;
\end{eqnarray*}
where $\eta$ is a small parameter.
\end{definition}
The horseshoes are going to be constructed on this Poincar\'e section. The 
auxiliary Poincar\'e section is defined differently from that in \cite{Li99}. 
\begin{definition}
The Poincar\'e section $\Sg_1$ is defined by the constraints:
\begin{eqnarray*}
& & z_1 =\eta ,\quad  -\eta < z_2 <\eta ,  \\
& & \sqrt{x^2+y^2} < \eta , \quad \|Q\| < \eta . 
\end{eqnarray*}
\end{definition}
The Poincar\'e map is defined as follows.
\begin{definition}
The Poincare map $P$ is defined as:
\[
P \ : \ U \subset \Sg_0 \mapsto \Sg_0, \quad P=P_1^0 \circ P_0^1,
\]
where
\[
P_0^1\ : \ U_0 \subset \Sg_0 \mapsto \Sg_1,\quad \forall \vQ \in U_0, 
\quad P_0^1(\vQ) = F^{t_0}(\vQ) \in \Sg_1,
\]
and $t_0=t_0(\vQ) >0$ is the smallest time $t$ such that 
$F^{t}(\vQ) \in \Sg_1$, and
\[
P_1^0\ : \ U_1 \subset \Sg_1 \mapsto \overline{\Sg_0} 
(=\Sg_0 \cup \pa \Sg_0),\quad \forall \vQ 
\in U_1, \quad P_1^0(\vQ) = F^{t_1}(\vQ) \in \overline{\Sg_0},
\]
and $t_1=t_1(\vQ) >0$ is the smallest time $t$ such that 
$F^{t}(\vQ) \in \overline{\Sg_0}$.
\label{dfp}
\end{definition}
The map $P_0^1$ has the explicit representation: Let $\vQ^0$ and $\vQ^1$
be the coordinates on $\Sg_0$ and $\Sg_1$ respectively, $y^0=0$, and 
$z_1^1=\eta$, then $t_0 = \frac{1}{\ga_1} \ln \frac{\eta}{z_1^0}$, and
\begin{eqnarray*}
x^1 &=& \bigg (\frac{z_1^0}{\eta}\bigg )^{\frac{a}{\ga_1}}x^0 \cos \bigg 
[ \frac{b}{\ga_1}\ln \frac{\eta}{z_1^0} \bigg ]\ , \\
y^1 &=& \bigg (\frac{z_1^0}{\eta}\bigg )^{\frac{a}{\ga_1}}x^0 \sin \bigg 
[ \frac{b}{\ga_1}\ln \frac{\eta}{z_1^0} \bigg ]\ , \\
z_2^1 &=& \bigg (\frac{\eta}{z_1^0}\bigg )^{\frac{\ga_2}{\ga_1}}z_2^0\ , \\
Q^1 &=& e^{t_0 L} Q^0\ .
\end{eqnarray*}
Let $\vQ_*^0$ and $\vQ_*^1$ be the intersection points of the homoclinic 
orbit $h_1$ with $\overline{\Sg_0}$ and $\Sg_1$ respectively. The discussion 
with respect to the other homoclinic orbit $h_2$ is the same. In a small neighborhood of $\vQ_*^1$, the map $P_1^0$ has an approximate representation. By virtue of the fact that the homoclinic orbit $h_1$ is a classical solution, the Well-Posedness Theorem \ref{WPT} implies that $F^t(\vQ_*^1)$ is $C^1$ in $t$. Thus for $\vQ^1$ in a small neighborhood of $\vQ_*^1$,
\begin{equation}
P_1^0(\vQ^1) = P_1^0(\vQ_*^1) + \L (\vQ^1 -\vQ_*^1) + 
\O (\| \vQ^1 -\vQ_*^1\|^2),
\label{app1}
\end{equation}
where
\[
\L (\vQ^1 -\vQ_*^1) = \pa_{\vQ} F^{t_1}(\vQ_*^1)\circ (\vQ^1 -\vQ_*^1) 
+\pa_t F^{t_1}(\vQ_*^1)\circ \frac{\pa t_1}{\pa \vQ^1}(\vQ_*^1)
\circ (\vQ^1 -\vQ_*^1) \ ,
\]
and $t_1=t_1(\vQ^1)$ is defined by the constraint that the $y$-coordinate 
of $F^{t_1}(\vQ^1)$ vanishes,
\[
F_y^{t_1}(\vQ^1) = 0 \ .
\]
Thus
\[
\pa_{\vQ} F_y^{t_1}(\vQ_*^1) + \pa_t F_y^{t_1}(\vQ_*^1) 
\frac{\pa t_1}{\pa \vQ^1}(\vQ_*^1) = 0 \ ,
\]
i.e.
\begin{equation}
\frac{\pa t_1}{\pa \vQ^1}(\vQ_*^1) = - \frac{1}{\pa_t F_y^{t_1}(\vQ_*^1)}
\pa_{\vQ} F_y^{t_1}(\vQ_*^1) \ .
\label{app2}
\end{equation}
Let $\vtQ^0 = \vQ^0 - \vQ^0_*$, and $\vtQ^1 = \vQ^1 - \vQ^1_*$, then 
$P_1^0$ has the approximate representation
\begin{equation}
\left ( \begin{array}{c} \tx^0 \\ \tz_1^0 \\ 
\tz_2^0 \\ \tQ^0 \\ \end{array} \right ) = \C 
\left ( \begin{array}{c} \tx^1 \\ \ty^1 \\ 
\tz_2^1 \\ \tQ^1 \\ \end{array} \right ) + \Xi ,
\label{rlp}
\end{equation}
where
\[
\Xi \sim  \O\bigg ((\tx^1)^2+(\ty^1)^2+(\tz_2^1)^2
+\|\tQ^1\|^2 \bigg ),
\]
\[
\C =  \left ( \begin{array}{cccc} c_{11} & c_{12}
& c_{13} & C_{14} \\  \\ c_{21} & c_{22}
& c_{23} & C_{24} \\  \\ c_{31} & c_{32}
& c_{33} & C_{34} \\  \\ C_{41} & C_{42}
& C_{43} & C_{44} \\ \end{array} \right ),
\]
in which $c_{jl}$ ($j,l=1,2,3$) are real constants, $C_{j4}$ ($j=1,2,3,4$) 
and $C_{4l}$ ($l=1,2,3$) are linear operators. 

\section{The Fixed Points of the Poincar\'e Map $P$}

As $t_0 \ra +\infty$, to the leading order, the fixed points of $P$ satisfy
\begin{equation}
\left ( \begin{array}{c} \hx^0 \\ 0 \\ 0 \\ \hQ^0 \\ \end{array} \right )
=\C \left ( \begin{array}{c} x^0_* \cos b t_0 \\ 
x^0_* \sin b t_0 \\ \hz^1_2 \\ 0 \\ \end{array} \right ), 
\label{ss1}
\end{equation}
where
\[
\hz^1_2 = e^{at_0} \tz^1_2, \quad \hx^0=e^{at_0}\tx^0, \quad  
\hQ^0=e^{at_0} \tQ^0 .
\]
Explicitly, the second and the third equations in (\ref{ss1}) are:
\begin{eqnarray}
x^0_*\bigg [ c_{21} \cos bt_0 + c_{22} \sin bt_0
\bigg ] + c_{23} \hz^1_2 =0, \nonumber \\
\label{ss2} \\
x^0_*\bigg [ c_{31} \cos bt_0 + c_{32} \sin bt_0
\bigg ] + c_{33} \hz^1_2 =0. \nonumber 
\end{eqnarray}
\begin{lemma}
$c_{23}$ and $c_{33}$ do not vanish simultaneously.
\end{lemma}

Proof. Notice that $W^u(Q_\e)$ is two-dimensional, and intersects $\Sg_0$
(or its extension to $-\eta < z_1 < \eta$) into a one-dimensional curve 
with tangent vector
\[
v= \C \left ( \begin{array}{c} 0 \\ 0 \\ 1 \\ 0 \\ \end{array} \right ).
\]
Notice also that for any $\vQ \in h_1$,
\begin{equation}
\mbox{dim} \{ \T_{\vQ}W^u(Q_\e) \cap \T_{\vQ}W^s(Q_\e) \} = 1, 
\label{dimc}
\end{equation}
where $\T_{\vQ}$ denotes the tangent space at $\vQ$. If $c_{23}$ and 
$c_{33}$ vanish simultaneously, then $v \in \T_{\vQ^0_*}W^s(Q_\e)$
which implies that
\[
\mbox{dim} \{ \T_{\vQ^0_*}W^u(Q_\e) \cap \T_{\vQ^0_*}W^s(Q_\e) \} = 2 
\]
which contradicts (\ref{dimc}). The lemma is proved. QED

Let 
\[
\Dl_1 = c_{21} c_{33}- c_{31} c_{23}, \quad
\Dl_2 = c_{22} c_{33}- c_{32} c_{23}.
\]
We assume that
\begin{itemize}
\item (A2). $\Dl_1$ and $\Dl_2$ do not vanish simultaneously.
\end{itemize}
Then (\ref{ss2}) has infinitely many solutions:
\begin{equation}
t^{(l)}_0 ={1 \over b}[l\pi -\varphi_1],\ \ l \in \mathbb{Z}; \label{ss4}
\end{equation}
where
\[
\varphi_1 = \arctan \{ \Dl_1 /\Dl_2 \}.
\]
Without loss of generality, we assume $c_{23}\neq 0$. 
Then, solving Eqs.(\ref{ss2}), we have
\begin{equation}
\hz^{(1,l)}_2 =  -x^0_* [c_{23}]^{-1} 
\{ c_{21} \cos b t^{(l)}_0 +c_{22} \sin b t^{(l)}_0 \}.
\label{ss5}
\end{equation}
Solving (\ref{ss1}), we have
\begin{eqnarray}
\hx^{(0,l)} &=& x^0_* \bigg [ c_{11} 
\cos b t^{(l)}_0 + c_{12}\sin b t^{(l)}_0 \bigg ] + c_{13}
\hz^{(1,l)}_2, \label{ss6} \\
\hQ^{(0,l)} &=& x^0_* \bigg [ C_{41} 
\cos b t^{(l)}_0 + C_{42}\sin b t^{(l)}_0 \bigg ] + C_{43}
\hz^{(1,l)}_2. \label{ss7}
\end{eqnarray}
Finally, by the implicit function theorem, there exist infinitely many 
fixed points of $P$, which have the approximate expressions given above 
\cite{Li99}. Specifically, we have 
\begin{theorem}
The Poincar\'e map $P$ has infinitely many fixed points labeled by $l$ 
($l \geq l_0$): 
\[
t_0 =t_{0,l} ,\ \hx^0 = \hx^0_l,\ \hQ^0 = \hQ^0_l,\ 
\hz^1_2 = \hz^1_{2,l},
\]
where as $l \ra +\infty$,
\begin{eqnarray*}
t_{0,l} &=& {1 \over b} [l\pi -\varphi_1] +o(1), \\
\hx^0_l &=& \hx^{(0,l)} + o(1), \\
\hQ^0_l &=& \hQ^{(0,l)} + o(1), \\
\hz^1_{2,l} &=& \hz^{(1,l)}_2 + o(1),
\end{eqnarray*}
in which $\hx^{(0,l)}$, $\hQ^{(0,l)}$ and $\hz^{(1,l)}_2$ are given in 
(\ref{ss6}),(\ref{ss7}),(\ref{ss5}).
\label{fixptthm}
\end{theorem}

\section{Existence of Chaos}

One can construct Smale horseshoes in the neighborhoods of the fixed points
of $P$.
\begin{definition}
For sufficiently large natural number $l$, we define slab $S_l$ in $\Sg_0$ 
as follows: 
\begin{eqnarray*}
S_l &\equiv& \bigg \{ \vQ \in \Sg_0 \ \bigg | \ \eta \exp \{ -\ga_1 
(t_{0,2(l+1)}-{\pi \over 2b})\} \leq \\
& & \tz^0_1(\vQ) \leq \eta \exp \{ -\ga_1 (t_{0,2l}-{\pi \over 2b})\}, \\
& & |\tx^0(\vQ)| \leq \eta \exp \{ -{1\over 2}a \ t_{0,2l} \}, \\
& & |\tz^1_2(P^1_0(\vQ))| \leq \eta \exp \{ -{1\over 2}a \ t_{0,2l} \}, \\
& & \|\tQ^1(P^1_0(\vQ))\| \leq \eta \exp \{ -{1\over 2}a \ t_{0,2l} \} 
\bigg \},
\end{eqnarray*}
where the notations $\tx^0(\vQ)$, $\tz^1_2(P^1_0(\vQ))$, etc. denote the 
$\tx^0$ coordinate of the point $\vQ$, the $\tz^1_2$ coordinate of the 
point $P^1_0(\vQ)$, etc..
\label{dfslab}
\end{definition}
$S_l$ is defined so that it includes two fixed points $p^+_l$ and $p^-_l$ 
of $P$ (Theorem \ref{fixptthm}). $S_l$, $P_0^1(S_l)$, and $\L P_0^1(S_l)$ 
are illustrated in Figure \ref{imag}, where $\L$ is defined in (\ref{app1}).
\begin{figure}
\vspace{3.0in}
\caption{An illustration of $S_l$, $P_0^1(S_l)$, and $\L P_0^1(S_l)$.}
\label{imag}
\end{figure}
\begin{figure}
\vspace{3.0in}
\caption{(a) shows one of the homoclinic 
orbits, and (b) shows the blow-up of the neighborhood 
of the saddle $Q_\e$.}
\label{A1}
\end{figure}

\nid
$\{ e_{\tx^0},e_{\tz^0_1},e_{\tz^0_2},{\bf e}_{\tQ^0}\}$ denotes the unit 
vectors along ($\tx^0,\tz^0_1,\tz^0_2,\tQ^0$)-directions in $\Sg_0$, 
$\{ e_{\tx^1},e_{\ty^1},e_{\tz^1_2},{\bf e}_{\tQ^1}\}$ denotes the unit 
vectors along ($\tx^1,\ty^1,\tz^1_2,\tQ^1$)-directions in $\Sg_1$, and 
under the linear map $\L$, $\{ e_{\tx^1},e_{\ty^1},e_{\tz^1_2},
{\bf e}_{\tQ^1}\}$ are mapped into $\{ \E_{\tx^1},\E_{\ty^1},\E_{\tz^1_2},
{\bf \E}_{\tQ^1}\}$. Let $e_\th$ be the unit angular vector at $P_0^1(p^+_l)$ 
of the polar coordinate frame on the ($\tx^1,\ty^1$)-plane. Let 
$E_\th = \L e_\th$, and we assume that
\begin{itemize}
\item (A3). $\ \mbox{Span} \bigg \{ e_{\tx^0},{\bf e}_{\tQ^0}, 
E_\th,\E_{\tz^1_2} \bigg \} = \Sg_0\ $.
\end{itemize}
Let $S_{l,\sg}= \sg \circ S_l$ where the symmetry $\sg$ is defined in 
(\ref{symm}). We need to define a larger slab $\hS_l$ such that 
$S_l \cup S_{l,\sg} \subset \hS_l$.
\begin{definition}
The larger slab $\hS_l$ is defined as
\begin{eqnarray*}
\hS_l &=& \bigg \{ \vQ \in \Sg_0 \ \bigg | \ \eta \exp \{ -\ga_1 
(t_{0,2(l+1)} -{\pi \over 2b})\} \leq \\
& & z^0_1(\vQ) \leq \eta \exp \{ -\ga_1 (t_{0,2l} -{\pi \over 2b})\},\\
& & |x^0(\vQ)-x^0_*| \leq \eta \exp \{ -{1\over 2} a\ t_{0,2l}\}, \\
& & |z^1_2(P^1_0(\vQ))| \leq |z^1_{2,*}| + \eta \exp \{ -{1\over 2} 
a\ t_{0,2l}\},\\
& & \|Q^1(P^1_0(\vQ))\| \leq \eta \exp \{ -{1\over 2}a\ t_{0,2l}\} \bigg \},
\end{eqnarray*}
where $z^1_{2,*}$ is the $z^1_2$-coordinate of $\vQ^1_*$. 
\end{definition}
\begin{definition}
In the coordinate system $\{ \tx^0,\tz^0_1,\tz^0_2,\tQ^0 \}$, the stable 
boundary of $\hS_l$, denoted by $\pa_s \hS_l$, is defined to be the boundary 
of $\hS_l$ along ($\tx^0,\tQ^0$)-directions, and the unstable 
boundary of $\hS_l$, denoted by $\pa_u \hS_l$, is defined to be the boundary 
of $\hS_l$ along ($\tz^0_1,\tz^0_2$)-directions. A stable slice $V$ in 
$\hS_l$ is a subset of $\hS_l$, defined as the region swept out through 
homeomorphically moving and deforming 
$\pa_s \hS_l$ in such a way that the part
\[
\pa_s \hS_l \cap \pa_u \hS_l
\]
of $\pa_s \hS_l$ only moves and deforms inside $\pa_u \hS_l$. The 
new boundary obtained through such moving and deforming of 
$\pa_s \hS_l$ is called the stable boundary of $V$, which is denoted by 
$\pa_s V$. The rest of the boundary of $V$ is called its unstable 
boundary, which is denoted by $\pa_u V$. An unstable slice of $\hS_l$,
denoted by $H$, is defined similarly.
\end{definition}
As shown in \cite{Li99}, under the assumption (A3), when $l$ is sufficiently 
large, $P(S_l)$ and $P(S_{l,\sg})$ intersect $\hS_l$ into four disjoint 
stable slices 
$\{ V_1,V_2 \}$ and $\{ V_{-1},V_{-2}\}$ in $\hS_l$. $V_j$'s 
($j=1,2,-1,-2$) do not intersect $\pa_s\hS_l$; moreover,
\begin{equation}
\pa_s V_i \subset P(\pa_sS_l), (i=1,2);\ 
\pa_s V_i \subset P(\pa_sS_{l,\sg}), (i=-1,-2).
\label{bdc}
\end{equation}
Let
\begin{equation}
H_j = P^{-1}(V_j),\ \ (j=1,2,-1,-2), \label{defhs}
\end{equation}
where and for the rest of this article, $P^{-1}$ denotes preimage of $P$. 
Then $H_j$ ($j=1,2,-1,-2$) are unstable slices. More importantly, the 
Conley-Moser conditions are satisfied as shown in \cite{Li99}. Specifically,
Conley-Moser conditions are:
\framebox[2in][l]{Conley-Moser condition (i):}
\[
\left \{ \begin{array}{l} V_j = P(H_j), \\ \pa_sV_j = P(\pa_sH_j), \ 
\ (j=1,2,-1,-2) \\ \pa_uV_j = P(\pa_uH_j).
\end{array}\right.
\]
\framebox[2in][l]{Conley-Moser condition (ii):} There exists a constant 
$0< \nu <1$, such that for any stable slice $V \subset V_j\ \ (j=1,2,-1,-2)$,
the diameter decay relation
\[
d(\tilde{V}) \leq \nu d(V)
\]
holds, where $d(\cdot)$ denotes the diameter \cite{Li99},
and $\tilde{V}=P(V\cap H_k), \ \ (k=1,2,-1,-2)$;
for any unstable slice $H \subset H_j\ \ (j=1,2,-1,-2)$, the diameter 
decay relation
\[
d(\tilde{H}) \leq \nu d(H)
\]
holds, where $\tilde{H}=P^{-1}(H\cap V_k), \ \ (k=1,2,-1,-2)$.

The Conley-Moser conditions are sufficient conditions for establishing
the topological conjugacy between the Poincare map $P$ restricted to
a Cantor set in $\Sg_0$, and the shift automorphism on symbols.

Let $\W$ be a set which consists of elements of the doubly infinite 
sequence form:
\[
a =(\cdot \cdot \cdot  a_{-2} a_{-1} a_0, a_1 a_2 \cdot \cdot \cdot ),
\]
where $a_k \in \{ 1, 2, -1, -2\}$; $k\in Z$. We introduce a topology in $\W$
by taking as neighborhood basis of 
\[
a^* =( \cdot \cdot \cdot a^*_{-2} a^*_{-1} a^*_0, a^*_1 a^*_2 
\cdot \cdot \cdot ),
\]
the set 
\[
W_j =  \bigg \{ a\in \W \ \bigg | \ a_k=a^*_k\ (|k|<j) \bigg \}
\]
\nid
for $j=1,2,\cdot \cdot \cdot $. This makes $\W$ a topological space.
The shift automorphism $\chi$ is defined on $\W$ by
\begin{eqnarray*}
\chi &:& \W \mapsto \W, \\
  & & \forall a \in \W,\ \chi(a) = b,\ \mbox{where}\ b_k=a_{k+1}.
\end{eqnarray*}
The shift automorphism $\chi$ exhibits {\em{sensitive dependence on 
initial conditions}}, which is a hallmark of {\em{chaos}}.

Let
\[
a =(\cdot \cdot \cdot  a_{-2} a_{-1} a_0, a_1 a_2 \cdot \cdot \cdot ),
\]
\nid
be any element of $\W$. Define inductively for $k \geq 2$ the stable
slices 
\begin{eqnarray*}
& & V_{a_0 a_{-1}}=P(H_{a_{-1}}) \cap H_{a_0}, \\
& & V_{a_0 a_{-1} ... a_{-k}}=P(V_{a_{-1} ... a_{-k}}) \cap H_{a_0} .
\end{eqnarray*}
\nid
By Conley-Moser condition (ii),
\[
d(V_{a_0 a_{-1} ... a_{-k}}) \leq \nu_1 d(V_{a_0 a_{-1} ... a_{-(k-1)}})
\leq ... \leq \nu_1^{k-1} d(V_{a_0 a_{-1}}).
\]
\nid
Then,
\[
V(a)=\bigcap^{\infty}_{k=1}V_{a_0 a_{-1} ... a_{-k}}
\]
\nid
defines a 2 dimensional continuous surface in $\Sg_0$; moreover,
\begin{equation}
\pa V(a) \subset \pa_u \hS_l. \label{inter1}
\end{equation}
\nid
Similarly, define inductively for $k \geq 1$ the unstable slices
\begin{eqnarray*}
& & H_{a_0 a_1}=P^{-1}(H_{a_1} \cap V_{a_0}), \\
& & H_{a_0 a_1 ... a_k}=P^{-1}(H_{a_1 ... a_k} \cap V_{a_0}) .
\end{eqnarray*}
\nid
By Conley-Moser condition (ii),
\[
d(H_{a_0 a_1 ... a_k}) \leq \nu_2 d(H_{a_0 a_1 ... a_{k-1}})
\leq ... \leq \nu_2^k d(H_{a_0}).
\]
\nid
Then,
\[
H(a)=\bigcap^{\infty}_{k=0}H_{a_0 a_1 ... a_k}
\]
\nid
defines a codimension 2 continuous surface in $\Sg_0$; moreover,
\begin{equation}
\pa H(a) \subset \pa_s \hS_l. \label{inter2}
\end{equation}
\nid
By (\ref{inter1};\ref{inter2}) and dimension count,
\[
V(a) \cap H(a) \neq \emptyset
\]
\nid
consists of points. Let 
\[
p\in V(a) \cap H(a)
\]
\nid
be any point in the intersection set. Now we define the mapping
\begin{eqnarray*}
\phi &:& \W \mapsto \hS_l, \\
  & & \phi(a) = p. 
\end{eqnarray*}
\nid
By the above construction,
\[
P(p)=\phi (\chi(a)).
\]
\nid
That is,
\[
P\circ \phi =\phi \circ \chi .
\]
\nid
Let 
\[
\La \equiv \phi (\W ),
\]
\nid
then $\La$ is a compact Cantor subset of $\hS_l$, and invariant
under the Poincare map $P$. 
Moreover, with the topology inherited from $\hS_l$ for $\La$,
$\phi$ is a homeomorphism from $\W$ to $\La$. Thus we have the theorem.
\begin{theorem}[Horseshoe Theorem]
Under the generic assumptions (A1)-(A3) for the 
perturbed nonlinear Schr\"odinger system (\ref{pnls}), there exists a 
compact Cantor subset 
$\La$ of $\hS_l$, $\La$ consists of points, and is invariant under $P$.
$P$ restricted to $\La$, is topologically conjugate to the shift 
automorphism $\chi$ on four symbols $1, 2, -1, -2$. That is, there exists
a homeomorphism
\[
\phi \ : \ \W \mapsto  \La,
\]
\nid
such that the following diagram commutes:
\begin{equation} 
\begin{array}{ccc}
\W &\maprightu{\phi} & \Lambda\\
\mapdownl{\chi} & & \mapdownr{P}\\
\W & \maprightd{\phi} & \Lambda
\end{array} 
\nonumber
\end{equation}
\label{horseshthm}
\end{theorem}

\section{Numerical Evidence for the Generic Assumptions}

\subsection{Generic Assumption (A1)}

Figure \ref{A1} shows a numerical result of Mark Winograd. It indicates 
that the homoclinic orbits are indeed tangent to the ($x,y$)-plane as 
$t \ra +\infty$. Specifically, Figure \ref{A1} (a) shows one of the homoclinic 
orbits, and Figure \ref{A1} (b) shows the blow-up of the neighborhood 
of the saddle $Q_\e$.

\subsection{Generic Assumptions (A2) and (A3)}
 
Numerical simulation of the generic assumptions (A2) and (A3) is planned 
for a future work.

\end{document}